\documentstyle[amsfonts,12pt]{article}
\newtheorem{thm}{Theorem}
\newtheorem{pro}[thm]{Problem}

\newtheorem{cor}[thm]{Corollary}
\newtheorem{lem}[thm]{Lemma}
\newtheorem{con}[thm]{Conjecture}

\def\endprf{\hfill\hbox{\kern1pt\vrule height6pt width4pt
depth1pt\kern1pt}\medskip}

\title{Zeros of Chromatic and Flow Polynomials of Graphs}
\author{Bill Jackson\\
Department of Mathematical and Computing Sciences,\\ Goldsmiths' College,
London SE14 6NW, England}
\date{20 April 2002}
\begin{document}
\maketitle

\begin{abstract}
We survey results and conjectures concerning the zero distribution of
chromatic and flow polynomials of graphs, and characteristic
polynomials of matroids.
\end{abstract}

\section{Introduction}
The study of
chromatic polynomials of graphs
was initiated by Birkhoff \cite{Bi} in 1912 and followed by Whitney
\cite{Wh1,Wh2} in 1932.
Inspired by the 4-Colour
Conjecture, Birkhoff and Lewis \cite{BL}, obtained results
concerning
the distribution of the
real zeros of
chromatic polynomials of planar graphs
and made the
stronger conjecture that
chromatic polynomials of planar graphs have no real zeros
greater than or equal to four.
Their hope was that results from analysis and algebra could be used to
prove their stronger conjecture and hence deduce that the 4-colour
conjecture was true. This has not yet occured: indeed the 4-colour
conjecture is now a theorem \cite{AH,AHK,RST}, but the stronger
conjecture of
Birkhoff and Lewis remains unsolved. Nevertheless, many beautiful
results have been obtained concerning the zero distibution of chromatic
polynomials both on  the real line  and in the complex plane,
and many other intriguing problems remain open.
I will  summarise these in Section \ref{chrom} below.

Flow polynomials were introduced by Tutte in \cite{Tutteflow}
as dual polynomials to chromatic polynomials. This duality holds only
for planar graphs, however,  and the zero distributions
 of chromatic and flow polynomials
for non-planar graphs seem to be quite different.
There are many similarities between
the zero distribution of flow polynomials for general graphs
and for planar graphs. This is not the case
for chromatic polynomials.
One possible reason for this is that chromatic polynomials of planar
graphs are strongly influenced by the fact that
the cycle space of a
planar graph has a basis consisting of a set of circuits which cover
every edge at most twice
(obtained by taking the boundaries of
all but one of the faces). The dual property, that
the cocycle space has a basis consisting of a set of cocircuits
which cover every edge at most twice holds for all graphs,
and can be obtained by taking the stars centred on all but one of the
vertices. Flow polynomials have received much less attention than
chromatic polynomials in the literature.
I will  summarise what little is known about their zero distribution
 in Section \ref{flow}.

Characteristic polynomials of matroids provide a common generalization
of chromatic and flow polynomials of graphs. I will describe how some
of the results from Sections \ref{chrom} and \ref{flow} can be extended to
this more general setting in Section \ref{matroid}.

All graphs considered are finite and may contain
loops and multiple edges.
We shall refer to graphs
without loops and multiple edges as {\bf simple graphs}.
We use $K_n$ and $K_{m,n}$ to denote the complete graph on $n$ vertices
and the complete bipartite graph with vertex sets of sizes $m$ and $n$,
respectively.

A graph $G$ is said to be {\bf non-separable} if $G$ is connected
and $G-v$ is connected for all $v\in V(G)$.
A subgraph $H$ of $G$
is a {\bf component}, respectively {\bf block},
of $G$ if it is a maximal connected, respectively non-separable,
subgraph of $G$.
The graph $G$ is $k$-connected if $|V(G)|\geq k+1$ and
$G-U$ is connected for all
$U\subseteq V(G)$ with $|U|<k$.
Given an edge $e$ in  $G$, we use $G-e$ to denote the graph
obtained from $G$ by deleting $e$, and $G/e$ to denote the graph
obtained from $G$ by contracting $e$. We say that $e$ is a {\bf bridge}
of $G$ if $G-e$ has more components than $G$.

\section{Chromatic Polynomials}\label{chrom}

Let $t$ be a positive integer and $G$ be a graph.
A {\bf proper $t$-colouring} of $G$ is an assignment of $t$
colours to
the vertices of $G$ in such a way that adjacent vertices of $G$ receive
different colours. If $G$ has a loop then $G$ does not have a proper
$t$-colouring for any $t$. If $G$ is loopless, then the
{\bf chromatic number}, $\chi(G)$ is the minimum
value of $t$  for which $G$ has a proper  $t$-colouring.
We use $P(G,t)$ to denote the {\bf chromatic polynomial} of $G$.
This is defined for integer values of $t\geq 1$ as the number of distinct
$t$-colourings of $G$. (Hence $P(G,t)\equiv 0$ if $G$ has a loop.)
The fact that $P(G,t)$ is a polynomial in $t$
and many other properties of $P(G,t)$ can be deduced from the following
elementary lemma.

\begin{lem}\label{fund1}
Let $G$ be a graph and $e$ be an edge of $G$. Then \[
P(G,t)=P(G-e,t)-P(G/e,t).\]
\end{lem}

Lemma \ref{fund1} also gives rise to a recursive procedure for
calculating chromatic polynomials. This calculation can often be
simplified using:

\begin{lem}\label{fund2}
Let $G$ be a graph and $G_1$ and $G_2$ be subgraphs of $G$ such that
$G_1\cup G_2=G$ and $G_1\cap G_2\cong K_r$. Then \[
P(G,t)=\frac{P(G_1,t)P(G_2,t)}{t(t-1)...(t-r+1)}.\]
\end{lem}
Note that it is NP-hard to determine the chromatic number
of a graph. Since the chromatic number of $G$ is one greater than
than the smallest positive integer root of $P(G,t)$, it is also NP-hard
to determine $P(G,t)$.

I will refer to the zeros of $P(G,t)$ as {\bf chromatic roots} of $G$.
The distribution of chromatic roots either on the real line
or in the complex plane is of
interest to
graph theorists principally because of its implication for the
integer chromatic roots and hence the chromatic number of a graph.
The distribution of chromatic roots  in the complex plane
is
also of independent interest to physicists because
of its relevence to phase transitions via the Potts model
partition function, see for example \cite[Section 1]{S1}.

\subsection{General Graphs}

Since  graphs can have arbitrarily large chromatic numbers
(for example $K_n$ has chromatic number $n$)
they can have arbitrarily large integer chromatic roots.
It follows,
from the definition of $P(G,t)$, however, that a graph $G$ can have no
integer chromatic root in the interval $[\chi(G),\infty)$.
One might hope that this result can be extended to real chromatic roots
but this is not the case. Woodall \cite{W1} has shown that complete
bipartite graphs (which have chromatic number two) have arbitrarily
large chromatic roots. More precisely, he showed that if $n$ is
large enough compared to $m$, then $K_{m,n}$ has chromatic roots
arbitarily close to $i$, for all integers $i$, $2\leq i\leq m/2$.

On the other hand, the following results due to
Tutte, Woodall and the author imply that the only
real numbers in $(-\infty,\frac{32}{27}]$ which are chromatic roots are
$0$ and $1$.

\begin{thm} \label{g1}
Let $G$ be a loopless graph with $n$ vertices,
$c$ components, and $b$ blocks which are not isolated vertices.
Then:
\begin{enumerate}
\item $P(G,t)$ is non-zero with sign $(-1)^{n}$ for
$t\in (-\infty,0)$ \cite{Tutte};
\item $P(G,t)$ has a zero of multiplicity $c$ at $t=0$ \cite{Tutte};
\item $P(G,t)$ is non-zero with sign $(-1)^{n+c}$ for
$t\in (0,1)$, \cite{Tutte};
\item $P(G,t)$ has a zero of multiplicity $b$ at $t=1$, \cite{W1};
\item $P(G,t)$ is non-zero with sign $(-1)^{n+c+b}$ for
$t\in (1,\frac{32}{27}]\approx (1,1.185]$, \cite{J};
\end{enumerate}
\end{thm}

The next result, due to  Thomassen and  Sokal, indicates that it is not
possible to extend Theorem \ref{g1}  beyond the
real interval $(-\infty, \frac{32}{27}]$ or into any region of the
complex plane.

\begin{thm}\label{dense1}
\begin{enumerate}
\item The real chromatic roots of graphs are dense everywhere in
$[\frac{32}{27},\infty)$ \cite{T1}
\item The complex chromatic roots of graphs are dense
everywhere in the complex plane.
  \cite{S2}.
\end{enumerate}
\end{thm}

\subsection{3-connected graphs}
We can construct 2-connected graphs with real chromatic roots in
$(1,2)$ as follows.
Let $G$ be a 2-connected graph
with an odd number of vertices.
It follows from Theorem \ref{g1}(c),(d) that
$P(G,1)=0$ and that $P(G,t)$ becomes negative when
$t$ increases from 1.
If, in addition, we assume that $G$ is bipartite, then $P(G,2)=2$
and hence $P(G,t)$ must have a zero in $(1,2)$. We may combine these
graphs with any other graph using Lemma  \ref{fund2}, with $r=2$,
to construct 2-connected non-bipartite graphs with chromatic roots
in  $(1,2)$.
It is conceivable, however, that Theorem \ref{g1} can be extended
if we add the hypothesis that $G$ is 3-connected.
(The same construction cannot be used to construct 3-connected
non-bipartite examples because $K_r$ is not bipartite when $r=3$.)

\begin{con}\label{3con}
Let $G$ be a loopless 3-connected graph.
\begin{enumerate}
\item $P(G,t)$ is non-zero with sign
$(-1)^{n}$ for $t\in (1,\alpha)$, where $\alpha\approx 1.781$ is the
chromatic root of $K_{3,4}$ in $(1,2)$.
\item If $G$ is not a bipartite graph with an odd number of vertices,
then $P(G,t)$
is non-zero with sign
$(-1)^{n}$ for $t\in (1,2)$.
\end{enumerate}
\end{con}
Conjecture \ref{3con}(b) is a slight strengthening of a
conjecture given in   \cite{J}.

The proof technique used in proving Theorem \ref{g1} is inductive,
using
Lemmas \ref{fund1} and \ref{fund2}. The main difficulty
in extending this technique to prove Conjecture \ref{3con}, is that
edge deletion and contraction may result in a graph
of connectivity 2, to which Lemma \ref{fund2} cannot be directly
applied.

There seems to be no hope of extending Theorem \ref{g1} into the complex
plane by adding a connectivity hypothesis since Sokal's result
Theorem \ref{dense1}(b) holds for chromatic roots of $k$-connected graphs for
any fixed $k$.

\subsection{Hamiltonian graphs}
Thomassen \cite{T0} suggested that the family of
hamiltonian graphs
may be another family for which Theorem \ref{g1} can be extended.
He made the following attractive conjecture:

\begin{con}\label{ham}\cite{T0}
If $G$ is hamiltonian and loopless then $P(G,t)$
is non-zero with sign
$(-1)^{n}$ for $t\in (1,2)$.
\end{con}

As evidence in favour of this conjecture, Thomassen
showed
in \cite{T2} that the zero free interval  of $(1,\frac{32}{27}]$
can be extended for graphs with a Hamilton path.

\begin{thm} \label{hampath}
Let $G$ be a loopless graph with $n$ vertices
and $b$ blocks.
If $G$ has a Hamilton path, then $P(G,t)$ is non-zero
with sign $(-1)^{n+b+1}$ for $t\in (1,\gamma]$ where
$\gamma\approx 1.296$ is the real root of $(2-t)^3=4(t-1)^2$.
\end{thm}
Thomassen \cite{T2} also constructed a sequence of graphs with
Hamilton paths
and with chromatic roots converging to $\gamma$.

We can use Lemmas \ref{fund1} and \ref{fund2} to show that
a smallest conterexample to Conjecture \ref{ham}
would be 3-connected, see \cite{T0}.
Thus Conjecture \ref{ham} would follow from Conjecture \ref{3con}(b).
Conjecture \ref{ham} would also follow from an affirmative
answer to the following intriguing
problem posed by Thomassen, using Lemmmas \ref{fund1} and
\ref{fund2}, see \cite{T0}.

\begin{pro}\label{hamred}\cite{T0}
Does every hamiltonian graph of minimum degree at least three
have an edge $e$ such that both $G-e$ and $G/e$ are hamiltonian?
\end{pro}

\subsection{Planar Graphs}\label{planar}
As mentioned above, chromatic polynomials were introduced by
Birkhoff and Lewis \cite{BL} in 1946 as a means of attacking the four
colour conjecture. Heawood \cite{H} had already proved that the chromatic number of any planar graph is at most five.
One of their main results in \cite{BL} is the
following nice generalization of Heawood's theorem.

\begin{thm}\label{bl1} Let $G$ be a loopless planar graph.
Then
$P(G,t)>0$  for
$t\in [5,\infty)$.
\end{thm}

Birkhoff and Lewis \cite{BL} conjecture that a similar
extension of the 4-colour theorem holds.

\begin{con}\label{bl2}
Let $G$ be a loopless planar graph.
Then $P(G,t)>0$ for $t\in [4,\infty)$.
\end{con}

One might hope that Theorem \ref{g1} could also be extended into the real
interval $(\frac{32}{27},4)$ for the special case of planar graphs.
The following result and conjecture of Thomassen \cite{T1}
indicate that this is probably not the case.

\begin{thm}\label{planerealdense}
The real chromatic roots of planar graphs are dense everywhere in
$[\frac{32}{27},3]$.
\end{thm}

\begin{con}\label{carst4}
The real chromatic roots of planar graphs are dense everywhere in
$[\frac{32}{27},4]$.
\end{con}

One may also hope that the complex plane may contain zero-free regions
for chromatic polynomials of planar graphs.
The following result of
Sokal \cite{S2} limits the location of any such region.

\begin{thm} \label{planecomdense}
The complex chromatic roots of planar graphs are dense
everywhere in the complex plane with the possible exception
of the disc $|t-1|< 1$.
\end{thm}

The planar graphs which Sokal uses to prove
Theorem \ref{planecomdense} are simply
two vertices joined by $p$ internally disjoint paths of length $d$.
He shows that as $p$ and $d$ vary, the complex chromatic roots of
graphs in this family are dense everywhere outside $|t-1|< 1$.
Graphs in the family also have roots inside the disc $|t-1|< 1$, but
it is still an open problem to determine whether the complex roots
of planar graphs are dense everywhere in the complex plane.
Theorem \ref{planecomdense} shows that Theorem \ref{bl1} cannot
be extended to the complex plane.

Much less is known about complex chromatic roots of planar
graphs of connectivity greater than 2.
Read and Tutte \cite{RT} have shown that the `bipyramids'
$C_n+\bar{K_2}$
are examples of 4-connected plane
triangulations with complex chromatic roots of unbounded modulus,
and Jacobsen, Salas, and Sokal \cite{JSS} have constructed
4-connected plane graphs which have complex chromatic roots with
real part greater than 4. Thus there is no obvious extension of
Conjecture \ref{bl2} to the complex plane, even for 4-connected
planar graphs.

\subsection{Plane Triangulations}\label{triangle}

A {\bf plane triangulation} is a loopless
plane graph in which all faces have size three.
They have a special significance in the study of chromatic roots
of planar graphs; for example,
we can use Lemma \ref{fund1} to show that Conjecture \ref{bl2}
is true if and only if it is true for plane triangulations.
We shall see below that the chromatic roots of plane
triangulations are much better behaved than those of planar graphs in
general.

Birkhoff and Lewis, and Woodall,  have
shown that Theorem \ref{g1}
can be extended beyond
$(-\infty,\frac{32}{27}]$ for plane triangulations.

\begin{thm}\label{triang1}
 Let $G$ be a  3-connected plane triangulation with $n$ vertices.
Then:
\begin{enumerate}
\item $P(G,t)$ is non-zero with sign $(-1)^{n}$ for
$t\in (1,2)$ \cite{BL};
\item $P(G,t)$ has a simple zero
at $t=2$ \cite{W1};
\item $P(G,t)$ is non-zero with sign $(-1)^{n+1}$ for
$t\in (2,\delta)$, where $\delta\approx 2.546$ is the chromatic
root of the octahedron in $(2,3)$ \cite{W2}.
\end{enumerate}
\end{thm}

One may construct infinite families of 3-connected plane triangulations
with a chromatic root at $\delta$ by combining the octahedron with
an arbitary plane triangultaion using Lemma \ref{fund2} with $r=3$.
In particular,
if we iteratively combine copies of the octahedron with itself,
we obtain an infinite family of plane triangulations whose only chromatic
root in $(2,3)$ is $\delta$.
On the other hand, the following beautiful conjectures of Woodall \cite{W3}
suggest that it may be possible
to extend Theorem \ref{triang1} for plane
triangulations of connectivity greater than three.

\begin{con} \label{woodcon1}
Let $G$ be a plane triangulation.
\begin{enumerate}
\item If $G$ is 4-connected then $P(G,t)$ has at most one zero in
$(2,\tau^2)$,
where $\tau=\frac{1+\sqrt{5}}{2}$ is the golden ratio, and
$\tau^2\approx 2.61803$.
\item If $G$ is 5-connected then $P(G,t)$ has
exactly one zero in
$(2,\theta)$, and
no zeros in
$(\theta,3)$, where $\theta\approx 2.61819$ is the chromatic
root of the icosahedron in $(2,3)$.
\end{enumerate}
\end{con}

\begin{con} \label{woodcon2}
\begin{enumerate}
\item Fot all $\epsilon>0$, there exist only finitely many
4-connected plane triangulations with a chromatic root in
$(2,\tau^2-\epsilon)$.
\item For all $\epsilon>0$, there exist only finitely many
5-connected plane triangulations with a chromatic root in
$(\tau^2+\epsilon,3)$.
\end{enumerate}
\end{con}

The following
intriguing result of Tutte \cite{Tutte2,Tutte3} shows that $\tau^2$ has
a special significance for chromatic polynomials of plane triangulations.

\begin{thm} Let $G$ be a plane triangulation with $n$ vertices. Then \\
$0\neq |P(G,\tau^2)|\leq \tau^{5-n}$.
\end{thm}

Note however that Tutte's result should not  be seen as evidence that
plane triangulations have chromatic roots close to $\tau^2$:
we have seen above that there exists an infinite family of
3-connected plane triangulations whose only chromatic root in $(2,3)$
is
$\delta\approx 2.546$.

The number $\tau^2$ is an element of an infinite sequence of real numbers,
called the {\bf Beraha numbers},
which seem to have a special significance for chromatic polynomials of
plane triangulations, see \cite{Tutte,BKW}.
For each integer $r\geq 2$,
let $b_r=2+2\cos\frac{2\pi}{r}$. Thus
$b_2=0, b_3=1, b_4=2, b_5=\tau^2, b_6=3$, and
$\lim_{r\rightarrow\infty}b_r=4$.

\begin{con} (Beraha \cite {B}, see \cite[14.6]{JT})
There exists a plane triangulation with a real chromatic root in
$(b_r-\epsilon,b_r+\epsilon)$ for all $r\geq 2$ and all $\epsilon>0$.
\end{con}

This conjecture would imply that there exist plane triangulations with
real chromatic roots arbitrarily close to 4.
Note that it is mistakenly stated in \cite{JT} that
Beraha and Kahane \cite{BK} have shown that this is true:
in fact they construct
plane triangulations with
complex chromatic roots arbitrarily close to 4.
Indeed it is an open problem to determine the supremum of the real
chromatic roots of plane triangulations. (Lemma \ref{fund1} can be used to
show that this will be equal to the supremum over all planar graphs, and
hence Conjecture \ref{carst4} would also imply that the supremum is 4.)
The largest real
chromatic root of a plane triangulation that I know of, 3.8267...,
is a root of a plane triangulation found by D.R. Woodall \cite{W3}.

One could also consider the supremum of the real chromatic roots of
planar bipartite graphs. Salas and Sokal conjecture that
this is $\tau^2$.

\begin{con}\cite{SS} Let $G$ be a planar bipartite graph. Then
$P(G,t)>0$ for $t\in [\tau^2,\infty)$.
\end{con}
They constructed
families of planar bipartite graphs
with chromatic roots tending to $\tau^2$ in \cite{SS}.

\subsection{Graphs of Bounded Degree}
Let $G$ be a connected graph and $\Delta(G)$ denote the maximum
degree of $G$.
Brooks \cite{Brooks} has shown that the chromatic number of $G$
is at most $\Delta(G)+1$ with equality if and only if $G$ is a
complete graph or an odd circuit.
Brenti, Royle and Wagner \cite{Brenti}
asked whether some form of Brooks's theorem could be true for the complex
chromatic roots of $G$.
More precisely, they asked whether there
exists a real function $f$  such that all complex chromatic roots of
$G$ lie in the disc $|t|<f(\Delta(G))$.
(This is, in fact, equivalent to an earlier conjecture for regular graphs
due to Biggs, Damerell and Sands \cite{Biggs}.)
The conjectures have recently been
verified by Sokal \cite{S1} using an exciting new proof technique.

\begin{thm}\label{Delta}
Let $G$ be a graph.
Then
$|P(G,t)|>0$ for all complex $t$ with $|t|\geq C \Delta(G)$,
where $C\approx 7.9639$.
\end{thm}

It is conceivable that Sokal's theorem is valid for some $C<2$,
although there are examples which show that we must take $C>1$,
see \cite[Section 7]{S1},
and hence Brooks's theorem cannot be extended to a disc in the complex
plane. It is possible however that the following extension of Brooks's
theorem is valid.
\begin{con}\cite{S3}
Let $G$ be a graph.
Then
$|P(G,t)|>0$ for all complex $t$ with $\Re(t)> \Delta(G)$.
\end{con}

A graph $G$ is said to be {\bf \boldmath $d$-degenerate}
if every subgraph of
$G$ has a vertex of degree at most $d$. Let $D(G)$ denote the
minimum value of $d$ for which $G$ is $d$-degenerate. Clearly
$D(G)\leq \Delta(G)$.
It is well known that $\chi(G)\leq D(G)+1$. It is not possible,
however, to extend this result on integer chromatic roots to the real
line, let alone the complex plane, since Thomassen \cite{T1} has
constructed 2-degenerate graphs with arbitrarily large real
chromatic roots.

Let $G$ and $H$ be graphs. We say that $H$ is a {\bf minor} of $G$
if $H$
can be obtained from $G$ by a sequence of edge deletions and contractions.
A {\bf simple minor} of $G$ is a minor which contains no loops or multiple edges.
If we modify the definition of $d$-degenerate by
replacing the condition  that
``every subgraph of $G$ has a vertex of degree at most $d$''
by the condition that
``every simple minor of $G$ has a vertex of degree at most $d$''
then Woodall \cite{W4} has shown that we do obtain a bound on the
size of the real chromatic roots of $G$.

\begin{thm}\label{graphminor} Let $G$ be a graph. Suppose that every
simple minor of $G$ has a vertex of degree at most $d$.
Then
$P(G,t)>0$ for all real $t\in (d,\infty)$.
\end{thm}

Since every simple planar graph has a vertex of degree at most five,
his result implies Theorem \ref{bl1} for all $t>5$.
Note that Theorem \ref{graphminor} cannot be extended
to a disc in the complex plane since even
Theorem \ref{bl1} does not extend to a disc in the complex plane.

The following conjecture of Sokal, which was inspired by
a remark of Shrock and Tsai \cite[p220]{ST}, suggests a possible way
to extend Theorem \ref{Delta}
to a larger family of graphs. Let {\boldmath $\Lambda(G)$}
denote the maximum number of edge-disjoint paths joining any pair of
vertices of $G$. Thus $\Lambda(G) \leq \Delta(G)$. Furthermore, it can
be seen that $D(G) \leq \Lambda(G)$ and hence $\chi(G)\leq \Lambda(G)+1$,
see \cite{JS}.

\begin{con}\label{Lambda}  \cite[Section 7]{S1}
There exists a constant $C$ such that
$|P(G,t)|>0$ for all complex $t$ with $|t|\geq C\Lambda(G)$.
\end{con}

\section{Flow Polynomials}\label{flow}

Let $\Gamma$ be an additive abelian group and $G$ be a graph.
Suppose we construct a digraph $\vec{G}$ by giving the edges of $G$ an
arbitrary orientation. For $U\subseteq V(G)$
and $\bar{U}=V(G)-U$, let $E^+(U)$ be the set of arcs from $U$ to
$\bar{U}$ in $\vec{G}$ and $E^-(U)= E^+(\bar{U})$.
Let $f:E(\vec{G})\rightarrow \Gamma$ and put
$f^+(U)=\sum_{e\in E^+(U)}f(e)$ and $f^-(U)=\sum_{e\in E^-(U)}f(e)$.
For $v\in V(G)$ let $f^+(v)= f^+(\{v\})$ and $f^-(v)= f^-(\{v\})$.
Then $f$ is a {\bf\boldmath $\Gamma$-flow}
for $G$, with respect to $\vec{G}$,
if $f^+(v)= f^-(v)$ for all $v\in V(G)$.
If in addition, $f(e)\neq 0$ for all $e\in E(G)$, then we say that
$f$ is a {\bf\boldmath nowhere-zero $\Gamma$-flow}
of $G$.
It can be seen that the condition $f^+(v)= f^-(v)$ for all $v\in V(G)$
is equivalent to the apparently stronger
condition that
$f^+(U)= f^-(U)$ for all $U\subseteq V(G)$. Thus, if $G$ has a
nowhere zero $\Gamma$-flow, then $G$ is bridgeless.
Since reversing
the orientation on an edge $e$ of $\vec{G}$ is equivalent to
replacing $f(e)$ by $-f(e)$,
the number of distinct nowhere-zero $\Gamma$-flows for
$G$ is independent of the chosen orientation $\vec{G}$ of $G$.

A {\bf\boldmath nowhere-zero $t$-flow}
of $G$ is a nowhere-zero ${\Bbb{Z}}$-flow, $f$,
such that $|f(e)|\leq t-1$
for all $e\in E(G)$. Tutte \cite{Tutteflow} has shown that
$G$ has a nowhere-zero $t$-flow if and only if $G$ has a
nowhere-zero ${\Bbb{Z}}_t$-flow. Furthermore, the number of distinct
nowhere-zero $\Gamma$-flows of $G$ is the same for all abelian
groups $\Gamma$ of the same order.
Note, however, that the number of
nowhere-zero ${\Bbb{Z}}_t$-flows of $G$ may differ from
the number of
nowhere-zero $t$-flows of $G$.
Nowhere-zero flows were introduced by Tutte \cite{Tutteflow} as a dual
concept to proper colourings. He showed that a connected
plane graph $G$
has a proper $t$-colouring if and only if its planar dual $G^*$ has a
nowhere-zero $t$-flow. The two concepts differ for non-planar graphs,
however. Indeed, whereas there exist loopless graphs which are not
$t$-colourable for arbitrarily large integers $t$,
the same is not true for bridgeless graphs and nowhere-zero $t$-flows.

\begin{thm}
Let $G$ be a bridgeless graph. Then:
\begin{enumerate}
\item $G$ has a nowhere-zero $6$-flow, (Seymour \cite{Seymour});
\item if $G$ has no 3-edge cuts then $G$ has a nowhere-zero $4$-flow,
(Jaeger \cite{Jaeger}).
\end{enumerate}
\end{thm}

Tutte  conjectures that both of these results can be
strengthened.

\begin{con} \label{tuttesflowcon}
Let $G$ be a bridgeless graph. Then
\begin{enumerate}
\item $G$ has a nowhere-zero $5$-flow \cite{Tutteflowpoly};
\item if $G$ has no 3-edge cuts then $G$ has a nowhere-zero $3$-flow.
\end{enumerate}
\end{con}
Conjecture \ref{tuttesflowcon}(b) was stated by Tutte in 1972.

Following Tutte \cite{Tutteflowpoly} we define the
{\bf\boldmath flow polynomial $F(G,t)$} of $G$ as the number of
distinct nowhere-zero ${\Bbb{Z}}_t$-flows of $G$ for any positive
integer $t$. Thus $F(G,t)\equiv 1$ if $E(G)=\emptyset$ and
$F(G,t)\equiv 0$ if $G$ has a bridge.
We shall refer to the zeros of
$F(G,t)$ as {\bf flow roots} of $G$.

By the above remarks, $F(G,t)$ is independent of the chosen
orientation
of $G$, and remains the same if we replace ${\Bbb{Z}}_t$ by any
other abelian group of order $t$. Note also that since the existence
of a nowhere-zero $t$-flow for $G$ implies the existence
of a nowhere-zero $(t+1)$-flow by definition, and is
equivalent to the existence of a nowhere-zero ${\Bbb{Z}}_t$-flow
as noted above,
we may deduce that if $F(G,t_0)\neq 0$ for some positive integer $t_0$,
then $F(G,t)\neq 0$ for all integers $t\geq t_0$.

We could also consider the polynomial $I(G,t)$ defined to be
the number of
distinct nowhere-zero $t$-flows of $G$. Kochol \cite{K} gives
relationships between $I(G,t)$ and $F(G,t)$, but these seem to
be the only results on $I(G,t)$ in the literature.
Attention has concentrated on $F(G,t)$ because  it
is dual to $P(G,t)$ for plane graphs.
Let $G$ be a connected plane graph and
$G^*$ be is its planar dual. There is a surjection from the
$t$-vertex-colourings of $G^*$ to the
nowhere-zero ${\Bbb{Z}}_t $-flows of $G$, such that each
nowhere-zero ${\Bbb{Z}}_t $-flow of $G$ has exactly $t$ pre-images,
see \cite{Tutteflowpoly}.
Thus $$F(G,t)=t^{-1}P(G^*,t).$$
We may use this identity to restate the results and conjectures
on chromatic roots of plane graphs in subsections \ref{planar}
and \ref{triangle} in terms of flow roots of plane graphs.
In the remainder of this section we shall survey results and
conjectures on flow roots of graphs which are not necessarily
planar.

We first state some fundamental recurrence relations for flow
polynomials.

\begin{lem}\label{flowfund1}
Let $G$ be a graph and $e$ be an edge of $G$.
\begin{enumerate}
\item If $e$ is a loop, then $F(G,t)=(t-1)F(G-e,t)$.
\item If $e$ is not a loop,
then $F(G,t)=F(G/e,t)-F(G-e,t)$.
\end{enumerate}
\end{lem}

\begin{lem}\label{flowfund2}
Let $G$ be a graph, $v$ be a vertex of $G$,
and $G_1$ and $G_2$ be edge-disjoint subgraphs of $G$ such that
$G_1\cup G_2=G$ and $G_1\cap G_2= \{v\}$. Then \[
F(G,t)=F(G_1,t)F(G_2,t).\]
\end{lem}

\begin{lem}\label{flowfund3}
Let $G$ be a 2-connected graph,
$v$ be a vertex of $G$,
$e$ be an edge of $G$,
and
$H_1$ and $H_2$ be edge-disjoint subgraphs of $G$ such that
$H_1\cup H_2=G-e$ and $H_1\cap H_2= \{v\}$.
Let $G_1$
be obtained from $G$ by contracting $E(H_{2})$, and define
$G_2$ analogously.
Then \[
F(G,t)=\frac{F(G_1,t)F(G_2,t)}{(t-1)}.\]
\end{lem}

\begin{lem}\label{flowfund4}
Let $G$ be a graph, $S$ be a 3-edge-cut of $G$,
and $H_1$ and $H_2$ be the components of $G-S$.
Let $G_1$
be obtained from $G$ by contracting $E(H_{2})$, and define
$G_2$ analogously.
Then \[
F(G,t)=\frac{F(G_1,t)F(G_2,t)}{(t-1)(t-2)}.\]
\end{lem}

\subsection{General Graphs}
Wakelin \cite{Wakelin} obtained the following analogue of
Theorem \ref{g1} for flow roots of graphs.

\begin{thm} \label{f1}
Let $G$ be a bridgeless graph with $n$ vertices,
$m$ edges, $b$ blocks, and no isolated vertices.
Then:
\begin{enumerate}
\item $F(G,t)$ is non-zero with sign $(-1)^{m-n+1}$ for
$t\in (-\infty,1)$;
\item $F(G,t)$ has a zero of multiplicity $b$ at $t=1$;
\item $F(G,t)$ is non-zero with sign $(-1)^{m-n+b+1}$ for
$t\in (1,\frac{32}{27}]$.
\end{enumerate}
\end{thm}

Using planar duality, we can restate Theorems \ref{planerealdense}
and \ref{planecomdense} in terms of flow
roots to deduce that
there is no obvious way
to extend Theorem \ref{f1} beyond the real
interval $(-\infty,\frac{32}{27}]$, or into the complex plane.

\begin{thm}
\begin{enumerate}
\item The real flow roots of planar graphs are dense everywhere in
$[\frac{32}{27},3]$.
\item The complex flow roots of planar graphs are dense
everywhere in the complex plane with the possible exception
of the disc $|t-1|< 1$.
\end{enumerate}
\end{thm}

We can construct 2-connected
graphs with real flow roots in $(1,2)$ as follows.
Let $G$ be a 2-connected graph
with $n$ vertices and $m$ edges and suppose that $m-n$ is odd.
It follows from Theorem \ref{f1}(b),(c) that
$F(G,1)=0$ and that $F(G,t)$ becomes negative when
$t$ increases from 1.
If, in addition, we assume that $G$ is Eulerian, then $F(G,2)=1$
and hence $F(G,t)$ must have a zero in $(1,2)$. We may combine these
graphs with any other graph using Lemma  \ref{flowfund3}
to construct non-Eulerian 2-connected
graphs with flow roots in  $(1,2)$.
It is conceivable, however, that Theorem \ref{f1} can be extended
if we add the hypothesis that $G$ is 3-connected.

\begin{con}
Let $G$ be a 3-connected graph with $n$ vertices,
and $m$ edges. Then:
\begin{enumerate}
\item $F(G,t)$ is non-zero with sign
$(-1)^{m-n}$ for $t\in (1,\phi)$, where $\phi\approx 1.749$ is the
flow root of $K_{5}$ in $(1,2)$.
\item
If $G$ is not an Eulerian graph with $m-n$ odd, then $F(G,t)$
is non-zero with sign
$(-1)^{m-n}$ for $t\in (1,2)$.
\end{enumerate}
\end{con}

Welsh \cite{Welsh} has conjectured that the dual form of
Conjecture \ref{bl2} holds for all graphs (with the
slight relaxation  that 4 is allowed to be a flow root).

\begin{con}\label{w}
Let $G$ be a bridgeless graph.
Then $F(G,t)>0$ for all real $t\in (4,\infty)$.
\end{con}

It is not even known whether this conjecture is true if 4 is replaced by
any large constant $C$. Note however that there seems to be no
obvious analogue of Conjecture \ref{w} for the complex plane:
by applying planar
duality to the remark given at the end of subsection \ref{planar},
we may deduce that there exist 3-connected
plane cubic graphs with complex flow roots of unbounded modulus,
and that there exist 3-connected plane graphs with
complex flow roots whose real part is greater than 4.

\subsection{Cubic Graphs}

A graph $G$ is {\bf cubic} if all its vertices have degree three.
They have a special significance in the study of flow roots
of graphs; for example,
we can use Lemma \ref{flowfund1} to show that Conjectures
\ref{tuttesflowcon}(a) and
\ref{w}
are true if and only if they are true for cubic graphs.
We shall see below that the flow roots of cubic
graphs are much better behaved than
those of graphs in
general.
Note that, using planar duality, we can restate the results and conjectures
on chromatic
roots of plane triangulations in terms of flow roots of planar
cubic graphs.
I have recently obtained the following extension of
(the dual form of) Theorem \ref{triang1} to cubic graphs which are not
necessarily planar.

\begin{thm}\cite{Jflow} \label{cubic1}
Let $G$ be a  3-connected cubic graph with $n$ vertices,
and $m$ edges.
Then
\begin{enumerate}
\item $F(G,t)$ is non-zero with sign $(-1)^{m-n}$ for
$t\in (1,2)$;
\item $F(G,t)$ has a zero of multiplicity $1$ at $t=2$;
\item $F(G,t)$ is non-zero with sign $(-1)^{m-n+1}$ for
$t\in (2,\delta)$, where $\delta\approx 2.546$ is the flow
root of the cube in $(2,3)$.
\end{enumerate}
\end{thm}

It seems likely that the
non-zero interval given in Theorem \ref{cubic1}(c)
may be extended beyond $\delta$ for cubic
graphs of `higher connectivity'. We use the following
concept to give a measure of `4-connectivity' in cubic graphs.
A graph $G$ is {\bf \boldmath cyclically $k$-connected}
if, whenever we can express $G$ as $G=G_1\cup G_2$, where
$E(G_1)\cap E(G_2)=\emptyset$ and $G_1$ and $G_2$ both contain
circuits, we must have $|V(G_1)\cap V(G_2)|\geq k$. (Cyclic
$k$-connectivity is the dual concept to
$k$-connectivity in
plane graphs.) Using this concept we may pose the following
strengthenings of Conjectures \ref{woodcon1}(a) and
\ref{woodcon2}(a).

\begin{con}
Let $G$ be a cyclically-4-connected cubic graph.
Then $F(G,t)$ has at most one
zero in $(2,\tau^2)$.
\end{con}

\begin{con}
For all $\epsilon>0$, there exist only finitely many
cyclically-4-connected cubic graphs with a flow root in
$(2,\tau^2-\epsilon)$.
\end{con}


\section{Characteristic Polynomials of Matroids}\label{matroid}

The reader may have wondered whether there is a more general
framework in which the duality between chromatic and flow polynomials
of plane graphs can be extended. Matroids
provide us with such a framework. I refer the reader to Oxley
\cite{OX} or Welsh \cite{W} for matroid definitions not explicitly
stated in this
article.

The {\bf \boldmath characteristic polynomial $C(M,t)$} of a
matroid $M=(E,r)$
is the polynomial in $t$ defined by
\[
C(M,t) = \sum_{A\subseteq E} (-1)^{\left\vert A \right\vert}
 t^{r(E)-r(A)}.\]

We can define the dual matroid $M^*$ for any matroid $M$.
We can also associate a pair of dual matroids to every graph $G$,
the cycle matroid $M_G$ and the cocycle matroid $M_G^*$.
Then
$C(M_G,t)=t^{-c}P(G,t)$, where $c$ is the number of components of $G$,
and
$C(M_G^*,t)=F(G,t)$.
Furthermore, when $G$ is a connected plane graph with planar dual
$G^*$, the cycle matroid of $G^*$ is equal to
the cocycle matroid of $G$. This gives the above mentioned
identity $F(G,t)=t^{-1}P(G^*,t)$.

The integer chromatic roots and flow roots of a graph both
occur as sequences of consecutive integers.
Examples are given  in \cite[p254]{W} to show
that this basic property may not
hold for the integer zeros of the chracteristic
polynomial of a matroid. Nevertheless, some of the above mentioned
results on real and complex chromatic roots do extend to matroids.
The following result from \cite{EHJ} gives a common
generalization of Theorems \ref{g1} and \ref{f1}.

\begin{thm} \label{m1}
Let $M$ be a loopless matroid with rank $r$ and $b$
components.
Then:
\begin{enumerate}
\item $C(M,t)$ is non-zero with sign $(-1)^{r}$ for
$t\in (-\infty,1)$;
\item $C(M,t)$ has a zero of multiplicity $b$ at $t=1$;
\item $C(M,t)$ is non-zero with sign $(-1)^{r+b}$ for
$t\in (1,\frac{32}{27}]$.
\end{enumerate}
\end{thm}

Let $M$ and $N$ be matroids. We say that $N$ is a {\bf minor} of $M$
if $N$
can be obtained from $M$ by a sequence of deletions and contractions.
A {\bf simple minor} of $M$ is a minor which contains no loops or
circuits of length two. Oxley \cite{OX1} has shown that if
every cocircuit of $M$ has size at most $d$ then
$C(M,t)>0$ for all real $t\in (d,\infty)$.
His proof, which uses induction on $|E(M)|$, can be used to prove
a stronger inductive statement which extends Theorem \ref{graphminor}.

\begin{thm} Let $M$ be a matroid. Suppose that every
simple minor of $M$ has a cocircuit of size at most $d$.
Then
$C(M,t)>0$ for all real $t\in (d,\infty)$.
\end{thm}

Applying this result to the special case of cographic matroids, we obtain:

\begin{cor} Let $G$ be a bridgeless graph. Suppose that every
3-edge-connected minor of $G$ has a circuit of length at most $d$.
Then
$F(G,t)>0$ for all real $t\in (d,\infty)$.
\end{cor}

Since every 3-connected graph $G$ has a circuit of length at
most $2\log_2|V(G)|$, this gives:

\begin{cor} Let $G$ be a bridgeless graph on $n$ vertices.
Then
$F(G,t)>0$ for all real $t\in (2\log_2 n,\infty)$.
\end{cor}

It is also possible that
Theorem \ref{Delta} and Conjecture \ref{Lambda}
can be extended to binary matroids.
Given a binary matroid $M$, let
$\Lambda(M)=\min_B\max_{C\in B}\{|C|\}$
where the minimum is taken over all bases $B$ of the cocycle space
of $M$. (It can be seen that if $M$ is the cycle matroid of a graph $G$
then $\Lambda(M)= \Lambda(G)$, see \cite{JS}.)

\begin{con} \cite{JS}
There exists a constant $D$ such that
for all loopless binary matroids $M$
and all complex numbers $t$ with $|t|\geq D\Lambda(M)$,
we have
$C(M,t)\neq 0$.
\end{con}

Applying this conjecture when $M$ is the
cographic matroid of a graph $G$ we obtain a
conjecture for flow roots of graphs.
Let $\Lambda^*(G)=\min_B\max_{C\in B}\{|C|\}$,
where the minimum is taken over all bases $B$ of the cycle space
of $G$.
\begin{con} \cite{JS}
There exists a constant $D$ such that
for all bridgeless graphs $G$ and all
complex numbers $t$ with $|t|\geq D\Lambda^*(G)$,
we have
$F(M,t)\neq 0$.
\end{con}

\subsection*{Acknowledgements}
I would like to thank Alan Sokal and Douglas Woodall for many helpful
remarks and inspiring conversations on chromatic roots which have greatly
contributed to this survey.

\end{document}